\documentclass{article}

\usepackage[english]{babel}
\usepackage[utf8]{inputenc}
\usepackage{amsmath}
\usepackage{amsthm}
\usepackage{amssymb}
\usepackage{graphicx}
\usepackage[colorinlistoftodos]{todonotes}
\usepackage{hyperref} 
\usepackage{color}
\usepackage{verbatim}

\theoremstyle{definition}

\newtheorem{puzzle}{Puzzle}

\title{Coins and Logic}
\author{Tanya Khovanova}
\date{\today}

\begin{document}

\maketitle

\textit{In memory of Raymond Smullyan.}

\begin{abstract}
We establish fun parallels between coin-weighing puzzles and knights-and-knaves puzzles.
\end{abstract}

\section{Knights and knaves. Genuine and fake coins.}

Raymond Smullyan wrote numerous logic books with many delightful logic puzzles. Usually the action happens on an island where two kinds of people live: \textit{knights} and \textit{knaves}. Knights always tell the truth, knaves always lie. Raymond Smullyan coined the names for knights and knaves in his 1978 book \textit{What Is the Name of This Book?} \cite{RSWhat}. Here is a puzzle from the book.
\begin{quote}\label{puzzle:kk}
\begin{puzzle} In this problem, there are only two people, A and B, each of
whom is either a knight or a knave. A makes the following
statement: ``At least one of us is a knave.'' 
What are A and B?
\end{puzzle}
\end{quote}

\textbf{Solution.} If A is a knave, then the statement is true, which makes it a contradiction. Therefore, A is a knight telling the truth and B must be a knave.

What does this have to do with coins? We can say that a fake coin is like a liar and a genuine coin is like a truth-teller. In all coin problems in this article, the balance scale has two pans. Consider the following coin puzzle. 
\begin{quote}
\begin{puzzle}\label{p:onefakelighter} You are given $N$ coins that look identical, but one of them is fake and is lighter than the other coins. All real coins weigh the same. You have a balance scale that you can use twice to find the fake coin. What is the greatest number $N$ of coins such that there exists a strategy that guarantees finding the fake coin in two weighings?
\end{puzzle}
\end{quote}

This puzzle is a slight variation of the first and the easiest coin puzzle that appeared in 1945 \cite{Schell}. Let us discuss some theory that will solve Puzzle~\ref{p:onefakelighter} and provide the framework for many coin problems.

\section{Some Coin Theory}

We start by looking at the outcomes of the weighings. The balance scale has two pans and the same number of coins is put on each pan to be weighed. The output of one weighing can be described as one of these three types:

\begin{itemize}
\item ``$=$''---when the pans are balanced,
\item``$<$''---when the left pan is lighter,
\item ``$>$''---when the right pan is lighter.
\end{itemize}

Suppose there is a strategy that finds a fake coin in two weighings. Suppose, for example, that the result of these two weighings is the outcome string $<=$ that points to a particular fake coin. There are nine possibilities for different outcomes: $==$, $=<$, $=>$, $<=$, $<<$, $<>$, $>=$, $><$, and $>>$. Each of them points to not more than one coin. Theoretically, it could be that some outcomes of two weighings do not point to any coin at all. In any case the number of coins can't be more than the number of outcomes. We proved that $N \leq 9$.

If we find the strategy for nine coins, the problem is solved. But first, let us introduce two types of strategies: \textit{adaptive} strategies in which each weighing can depend on the results of all previous weighings, and \textit{oblivious} (or non-adaptive) strategies in which all the weighings must be specified in advance. 

Now we introduce itineraries. During any weighing, a coin's presence on the \textbf{l}eft pan is denoted by L; a coin's presence on the \textbf{r}ight pan is denoted by R; and a coin not participating (one that is left \textbf{o}utside of the weighing) is denoted by O. After all the weighings, every coin's path can be described as a string of Ls, Rs, and Os. Such a string is called the coin's \textit{itinerary}.

In an oblivious strategy the itinerary of every coin is known in advance. As we noted with outcomes, two coins can't have the same itineraries. Let us prove this by contradiction. Suppose two coins have the same itineraries. That means they are together in the same group in every weighing. It follows that if one of them is fake, we can't identify which one of them it is. There are nine possible different itineraries. Therefore, we can't have an oblivious strategy that finds the fake coin in two weighings if there are more than nine coins. But we already knew that.

Itineraries in an adaptive strategy are trickier. Suppose in the first weighing we compare three coins versus three coins and the weighing balances. We can find the fake coin in the second weighing if we compare two coins that were not on the scale in the first weighing. That means the three coins that were on the left pan during the first weighing have the same itineraries. The same is true for the three coins that were on the right pan during the first weighing. We found the fake coin, but the itineraries are not unique. That's the bad news. The good news is that something is unique. The itinerary of the fake coin can't be shared in this strategy with any other coin. This is for the same reason as before: the fake coin has to be distinguishable from other coins. We call the itinerary of the fake coin the \textit{self-itinerary}. 

Given an adaptive strategy every coin has a unique self-itinerary that finds this particular coin if it is fake. Moreover, distinct coins have distinct self-itineraries. We can prove this by contradiction. If two coins have the same self-itineraries, then any adaptive strategy would proceed the same way whether the first or the second coin is fake. That means these two coins will always be together and the strategy can't distinguish which of them is fake. 

If we go back to an oblivious strategy, where the itineraries are decided in advance, the self-itinerary is the same as the itinerary.

Let us go back to our sample outcome $<=$, which points to a coin. The fake coin that corresponds to this outcome has to have the self-itinerary $LO$. And vice versa, given a self-itinerary we can find the outcome string that points to this coin. Thus, if we have nine coins, we have a one-to-one correspondence between outcomes and self-itineraries: $<<$, $<>$, $<=$, $><$, $>>$, $>=$, $=<$, $=>$, $==$ point to LL, LR, LO, RL, RR, RO, OL, OR, OO correspondingly.

Now the magic happens. If we are looking for an oblivious strategy we must have these itineraries. To be precise we assign itineraries LL, LR, LO, RL, RR, RO, OL, OR, OO to coins numbered 1 through 9 correspondingly. The first weighing has coins 1, 2, and 3 on the left pan and coins 4, 5, and 6 on the right pan. The second weighing has coins 1, 4, and 7 on the left pan, and coins 2, 5, and 8 on the right pan. Hooray! It works. 

Where is the magic? The magic is in the fact that the itineraries do not contradict the rule of the balance scale: the same number of coins on each pan in each weighing. We produced an oblivious strategy with three coins on each pan in each weighing. The magic is not very surprising. Our itineraries are symmetric with respect to switching L and R. That means the pans should have the same number of coins. The magic is in mathematics!

Anyway, we proved that $N$ can't be more than 9, and found an oblivious strategy for 9 coins. Problem solved.

\section{Normals and Chameleons}

Now we get to normals. In logic puzzles \textit{normals} are not obsessive with telling a lie or the truth. Like normal people in life, they sometimes lie and sometimes tell the truth. Here is another puzzle from Raymond Smullyan's book \cite{RSWhat} with the wonderful self-referencing title:

\begin{quote}
\begin{puzzle} Consider a married couple, Mr.~and Mrs.~A. It is known that either both of them are normal, or one of them is a knight and the other a knave. They make the following statements:
Mr.~A: My wife is not normal.
Mrs.~A: My husband is not normal.

What are Mr.~and Mrs.~A?
\end{puzzle}
\end{quote}

\textbf{Solution.} Suppose Mr.~A is a knave. It follows from his lie that his wife is normal, which is a contradiction. Similarly, Mrs.~A can't be a knave. It follows that they are both normal.

\begin{figure}
\begin{center}
  \includegraphics[scale=0.4]{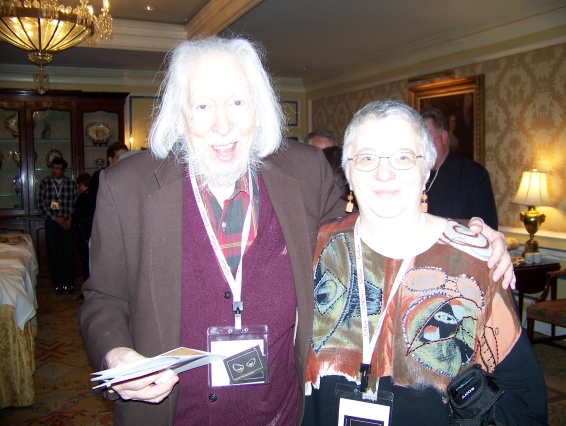}
\end{center}
\label{fig:RSandTK}
\caption{Raymond Smullyan and Tanya Khovanova at G4G8 in 2008}
\end{figure}

I was brought up on Raymond Smullyan’s books. I was so happy when I first met him at the Gathering for Gardner 8 in 2008, where he played a kissing trick on me \cite{TKRSblog}. Figure~\ref{fig:RSandTK} is a picture of us from the G4G8 conference.

Back to normals. The following puzzle will be useful later.

\begin{quote}
\begin{puzzle} What one statement can a normal make to prove their identity?
\end{puzzle}
\end{quote}

\textbf{Answer.} I am a liar.

What happens if a normal doesn't want anyone to know who they are? Then the interrogator is in trouble. The normal person can consistently behave as either a truth-teller or a liar.

What is a coin analogue of a normal person? 

The \textit{chameleon} coin can pretend to be real or fake (lighter). The chameleon was introduced and studied in my paper with Konstantin Knop and Oleg Polubasov \cite{KKP}. 

Did you notice that coins and people have different names in these puzzles. It is conceivable to hear someone call genuine coins truth-tellers and fake coins liars. But, really, it would be strange to call chameleon coins, normal coins. There is nothing normal about them.

What is a good question to ask about a chameleon? Suppose there are $N$ coins and one of them is a chameleon, while the others are genuine. Can we find the chameleon? Actually, no. If the chameleon always pretends to be real then it can't be found. Likewise, the chameleon can hide among fake coins never to be found. This is similar to a normal person who can pretend to be either a truth-teller or a liar. Luckily in most logic puzzle there are several people, and the normal can be ratted out by someone else. 

Suppose we have a mixture of real coins and fake coins in addition to one chameleon. We can't find the chameleon. But can we find the fake coin in the presence of a chameleon? No, again. If the chameleon pretends to be fake, we can't prove that a particular coin is fake. What do we do? Here is the simplest possible question we can ask.

\textbf{Question.} We have $N$ total coins. All but two of them are genuine and weigh the same, but two of them are not genuine. One is a classic fake that is lighter, and one is a chameleon. The goal is to find two coins that are guaranteed to include the classic fake in the smallest number of weighings.

If the chameleon always pretends to be fake, the two coins that we will find at the end will be the chameleon and the fake. If the chameleon sometimes pretends to be real, we might be able to find the fake coin, which is better than our goal.

Our paper \cite{KKP} is devoted to this question. Here is the summary. There is a bound for how fast we could find our two coins: it can't be faster than finding two fake coins among $N$ coins. This is because the chameleon can pretend to be fake. The starting point for many coin weighing problems is to count the number of outcomes. If there are $w$ weighings, then there are $3^w$ possible outcomes. In the problem of two fake coins each outcome points to two coins. This means that the total number $N$ of coins that can be processed is bounded by  
\[\binom{N}{2} \leq 3^w.\]
This bound is often called the information-theoretic bound (ITB). As you will see in a moment this bound is very good: it is close to the real answer.

The computational summary is in Table~\ref{tbl:itbound}. The table shows the results for up to ten weighings. The second row is the ITB we just calculated. The third row shows the results for two fake coins from \cite{KP} and the fourth row for one fake and one chameleon from \cite{KKP}. The numbers in bold are exact, derived from an exhaustive computer search. The other numbers in the last two rows represent current best algorithms. As you can see the ITB for two fake coins is very close to the best known algorithm. You can also see that the chameleon makes life more difficult.

\begin{table}[h!]
  \begin{center}
\begin{tabular}{| r | r | r | r | r | r | r | r | r | r | r |}
  \hline                       
  \# weigh-s & 1 & 2 & 3 & 4 & 5 & 6 & 7 & 8 & 9 & 10\\
  ITB & 3 & 4 & 7 & 13 & 22 & 38 & 66 & 115 & 198 & 344\\
  FF alg. & \textbf{3} & \textbf{4} & \textbf{7} & \textbf{13} & \textbf{22} & \textbf{38} & 65 & 113 & 194 & 341 \\
  cham. alg.  & \textbf{2} & \textbf{4} & \textbf{6} & \textbf{11} & \textbf{20} & 36 & 60 & 108 & 180 & 324\\
\hline
\end{tabular}
  \end{center}
\label{tbl:itbound}
\caption{ITB and results for a small number of weighings}
\end{table}

\section{Alternators}

Normal people are unpredictable. They can lie or tell the truth at any time. We can simplify normals by making them more deterministic. \textit{Alternators} switch between telling the truth and lying. Here is an alternator puzzle, whose provenance is not known.

\begin{quote}
\begin{puzzle} Folks living in Trueton always tell the truth. Those who live in Lieberg, always lie. People living in Alterborough alternate strictly between truth and lie. One night 911 received a call: ``Fire, help!'' The operator couldn't identify the phone number, so he asked, ``Where are you calling from?'' The reply was Lieberg.
Assuming no one had overnight guests from another town, where should the firemen go?
\end{puzzle}
\end{quote}

There is a similarity to the previous puzzle, where normals identify themselves in one statement. In this case, only an alternator can say that they are from Lieberg. The solution to the puzzle doesn't end here. To finish it, we need to say that because the last statement, which identified the town, was a lie, the first statement must have been true. That means there is indeed a fire, and help needs to be dispatched.

Here is another puzzle related to alternators.

\begin{quote}
\begin{puzzle} How can we check in two questions if a person is a knight, a knave or an alternator?
\end{puzzle}
\end{quote}

There are probably many ways. A standard approach is to ask if two plus two is four, twice. In particular, a lone alternator can always be found, unlike a normal person.

Logic puzzles with alternators are often easier than with normals, especially if there is an alternator that makes many statements. What is a coin analogue of the alternator? Of course, it is a coin that switches its weight between being real and being lighter after each weighing. Such a coin is called, unimaginatively,  an \textit{alternator}. This is the first time the names for coins and people match.

I do not know the history of people alternators. It is likely that in logic the normals appeared before the alternators. It is not surprising that coin puzzles followed suit: The alternator coins appeared after the chameleons. 

Anyway, I am running a math club for middle school students at MIT called PRIMES STEP. This is a junior program that branched out from another MIT youth program called PRIMES. Surprisingly, PRIMES is an abbreviation: Program for Research in Mathematics, Engineering, and Science for High School Students. A posteriori, it is less surprising that STEP is also an abbreviation: Solve–Theorize–Explore–Prove. I am the Head Mentor for both programs. 

In the 2015-2016 academic year, I gave the alternator puzzle as a research problem to my PRIMES STEP students who were in middle school: grades 6-8.

My hope was that they could find low-hanging fruits: check the small number of coins and find some trivial bounds. The students did much more.

Unlike the chameleon coin, one alternator coin among the real coins can always be found. We can use a strategy of finding a fake coin among real coins and repeat every weighing twice. Remember that human alternators also reveal themselves with the double strategy.

Let us have another logic puzzle. 

\begin{quote}
\begin{puzzle} It is known that A is an alternator. A makes the following two statements about B: B is a knave; B is an alternator. What is B?
\end{puzzle}
\end{quote}

\textbf{Solution.} This is not a kosher puzzle. We cannot deduce what B is. It depends on the starting state of A. If A starts with the truth, then B is a knave. If A starts with a lie, then B is an alternator.

Let us try again:

\begin{quote}
\begin{puzzle} A makes the following two statements about B: B is a knave; B is an alternator. To which B replies: A is an alternator. What are A and B?
\end{puzzle}
\end{quote}

\textbf{Solution.} A can't be a truth-teller as the two statements contradict each other. If A is a liar, then B is a truth-teller, leading to a contradiction. Therefore, A is an alternator, which means that B's statement is true. B can't be a liar. Therefore, B is an alternator.

As with the alternator people, it is useful to know the starting state of the alternator coin: whether the first time the coin is on the scale it pretends to be fake or real. Mathematically, we have three cases to research: 

\begin{enumerate}
\item the starting state is known and lighter, 
\item the starting state is known and real, 
\item the starting state is unknown.
\end{enumerate}

As usual, it is useful to count the number of possible outcomes of $n$ weighings. Didn't we already conclude that this number is $3^n$? Yes, we did. But in the case of alternators not every outcome string is possible: it can't have two imbalances in a row. Indeed, after an imbalance, the alternator changes weight and pretends to be real. The next weighing has to balance, no matter whether the alternator is on the scale or not. 

The question of counting outcomes that do not have two imbalances in a row is a combinatorial question, and the answer is supplied by the coolest sequence my students never heard of: the Jacobsthal numbers \cite{OEIS}. Here is the sequence:

\[0, \quad 1, \quad  1, \quad  3, \quad  5, \quad  11, \quad  21, \quad  43,  \quad 85,  \quad 171,  \quad \ldots\]

The sequence is defined by the initial terms and a recursion:

\begin{itemize}
\item $J_0 = 0$, $J_1 = 1$,
\item $J_n = J_{n-1} + 2 J_{n-2}$, for $n > 1$.
\end{itemize}

I leave it to the reader to see that the recursion must be this one. Alternatively, check the proof in the paper I wrote with the PRIMES STEP students \cite{STEPAlternators}.

Actually I wasn't precise. The number of outcomes depends on the starting state. If the alternator starts in the real state, the first weighing can't be an imbalance. So the count of $n$ outcomes is: 

\begin{itemize}
\item lighter state: $J_{n+2}$,
\item real state: $J_{n+1}$.
\end{itemize}

The cool part is that there is an adaptive strategy that processes this many coins. The number of outcomes is the precise bound.

What about the unknown state? The number of outcomes should be the same as the outcomes that are good for both the lighter state and the real state. Actually any outcome for the lighter state can also work for the unknown state. So the number of possible outcomes of $n$ weighings in an unknown state is $J_{n+2}$: the same as for the lighter state. Does this mean that we can process $J_{n+2}$ coins in $n$ weighing when the alternator starts in the unknown state? This can't be true. The unknown state can't be easier than the real state. What are we missing here? When the alternator coin starts in the unknown state, if it is ever on the scale, we will not only find the coin, but after we find it, we will know the initial starting state. Roughly speaking, each coin matches two different outcomes. Roughly, because it might be possible that the coin is never on the scale, and we find it by proving that all the other coins are real. This way we do not know the starting state of the alternator. The good news---there can only be one such coin. That means we can't process more than $(J_{n+2}+1)/2$ coins in the unknown state. 

On the other hand, the number of coins in the unknown state that we can process can't be more than the number of coins in the real state, that is $J_{n+1}$. It follows that we can't process more unknown coins than the minimum of these two numbers $(J_{n+2}+1)/2$ or $J_{n+1}$.  The cool property of the Jacobsthal numbers is 
\[J_{n+1}=2J_{n}+(-1)^{n}\,\]
that is, the next Jacobsthal number is almost twice the previous number. In any case the minimum is $J_{n+1}$ and this is how many unknown coins we can process in an adaptive strategy. You can find more details in our paper \cite{STEPAlternators}.

\begin{figure}
\begin{center}
  \includegraphics[scale=0.2]{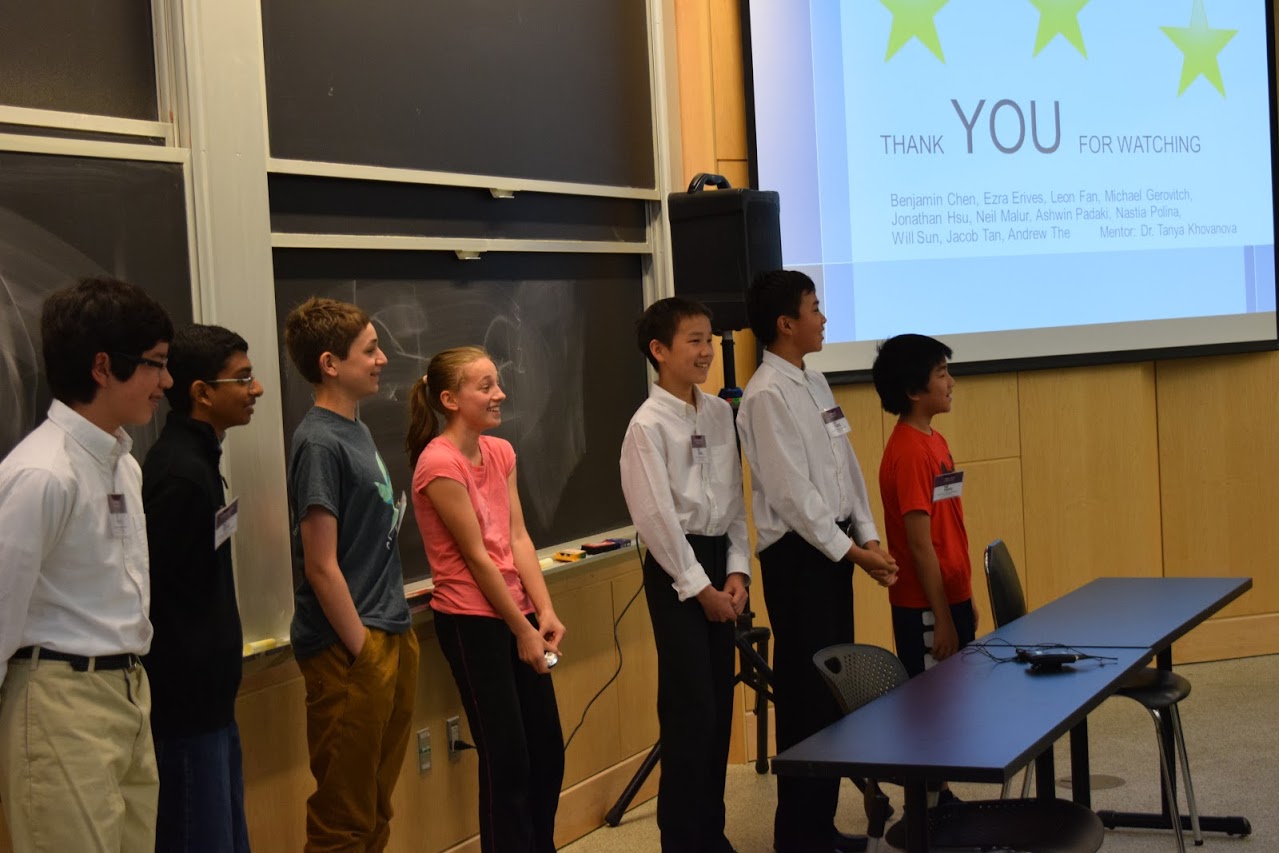}
\end{center}
\label{fig:PS}
\caption{PRIMES STEP students give a talk.}
\end{figure}

In Figure~\ref{fig:PS} the PRIMES STEP students give a talk about the research at the 2016 PRIMES conference.

\section{Alternator: an oblivious strategy}

What about an oblivious strategy for the alternator? The STEP students didn't find it. Is it possible to do? After the program was over I sent our PRIMES STEP paper to my coin co-author Konstantin Knop and we started to discuss the oblivious strategy. One of the problems is that the outcomes do not uniquely define self-itineraries. For example, consider the outcome $<==<$. The alternator has to be on the scale three times: in the first and the last weighing and exactly one more time in between, after which it needs to change its state from real to light. There are four possible itineraries: LLOL, LROL, LOLL, and LORL. Which one should we choose?

We can try to choose itineraries in a systematic way. For example, we can choose LLOL. We can say that after the coin changes weight to real we put it on the same scale again in the next weighing. Using this technique we found oblivious strategies for known starting states that matched in power the adaptive strategies, but the unknown state was more complicated \cite{KK}.

Consider the outcome $\neq = \neq = \neq$, where $\neq$ marks an imbalance. The coin matching this outcome must start in the light state and must be on the scale for every weighing. The same coin in the real state must have an outcome: $= \neq = \neq =$. That means eight outcomes of type $\neq = \neq = \neq$ should correspond to the same coins as the four outcomes of type $= \neq = \neq =$. It follows that we have outcomes that can't be matched to coins. The oblivious strategy for the unknown starting state that is as good as an adaptive strategy doesn't exist. What is the best oblivious strategy? This is still an open problem. 
See \cite{KK} for more details.

How do adaptive and oblivious strategies transfer to logic? For an oblivious strategy you prepare your questions in advance. Let's try a puzzle.

\begin{quote}
\begin{puzzle} You are on an island where knights, knaves and alternators live. You meet two people. You want to figure out who they are in the smallest number of questions. The only type of question you can ask to one of them is: Is the other person a knight? You can replace “knight” in this question with either “knave” or “alternator.”

What is your oblivious strategy?
\end{puzzle}
\end{quote}

\textbf{Hint.} There are nine possibilities for who those two people can be. That means you need at least four questions. It is possible to prepare four questions in advance to figure out who is who. How? I leave it to you.

\section{Heavier versus Lighter}

Traditionally fake coins are lighter. The reason is historical: cheaper metals were lighter.

We mathematicians like symmetry. We can imagine a fake coin that is heavier than the real one. Finding such a coin is not an interesting problem. We can use the same strategy as for a lighter fake coin, we just need to interpret the outcomes in the opposite way: $<$ means the fake heavier coin is on the right pan. As a result most math problems about counterfeit coins assume that the fake coin is lighter.

There are, on the other hand, puzzles where you need to find a fake coin that is heavier or lighter. The most famous is the 12-coin problem \cite{Dyson} first introduced in 1945.

\begin{quote}\label{puzzle:12coins}
\begin{puzzle} There are 12 identical-looking coins. One of them is fake: it is either lighter or heavier than a real coin. Find the fake coin and tell if it is lighter or heavier by using a balance scale three times. 
\end{puzzle}
\end{quote}

This puzzle is so famous that it shows up on top of searches for coin problems. There are so many papers and books explaining the solution that I don’t have to go into detail. I will just provide an oblivious strategy.

\textbf{Solution.} First weighing: coins 1, 2, 3, 4 versus 5, 6, 7, 8. Second weighing: coins 1, 8, 9, 11 versus 3, 4, 5, 10. Third weighing: coins 3, 6, 10, 12 versus 1, 4, 6, 11. In the language we are discussing we should translate this strategy to itineraries: LLR, LOO, LRL, LRR, RRO, ROR, ROL, RLO, OLL, ORO, OLR, OOL for coins 1 through 12 correspondingly.

The important thing to notice here is that there are no two coins that are opposite each other each time they are on the scale. There is a deep reason why this should be the case. If there were two coins like that and one of them was fake, we wouldn't be able to differentiate between the two. The strategy will lead to the conclusion that either one of them is fake and lighter or the other one is fake and heavier. 

There is a variation of this puzzle where you just need to find the fake coin.

\begin{quote}\label{puzzle:13coins}
\begin{puzzle} There are 13 identical-looking coins. One of them is fake: it is either lighter or heavier than a real coin. Find the fake coin by using the balance scale three times. 
\end{puzzle}
\end{quote}

The strategy is similar to the one in the previous problem, but we just need to add an itinerary for the 13th coin, which is OOO. If all the weighings balance, then the last coin is fake and we do not know whether it is heavier or lighter.

Before moving to logic, I should mention that we can expand this problem to more weighings. If we need to find and identify the fake coin, we can process $(3^w-3)/2$ coins in $w$ weighings, and one more coin if we just need to find it.

If we go back to logic, what would be the analogues of lighter versus heavier fake coins? They both lie.

It might look like logic is binary: truth versus lie, while, on the other hand, coin weighing is ternary. If you think about it, then in real life lies can take many forms. For example, they can say that two plus two is three, or two plus two is seventeen. 

For our logic puzzles, let’s introduce two different types of knaves.

 Suppose an \textit{exaggerating liar} increases any number he mentions by 1, and a \textit{diminishing liar} decreases any number by 1.

\begin{quote}
\begin{puzzle} On an island there are three types of people: exaggerating and diminishing liars and truth-tellers. You meet three islanders A, B, and C. They say the following:
A: There are three liars among us.
B: There is one liar among us.

Who is who?
\end{puzzle}
\end{quote}

\textbf{Solution.} A's statement must be a lie. Therefore, A can only be an exaggerating liar. Hence, there are two liars in the group. B must be a diminishing liar and C a truth-teller.

Let us consider another coin problem from the book 
\textit{Mathematical Circles: Russian Experience} \cite{FGI}: 

\begin{quote}
\begin{puzzle}\label{puzzle:h-or-l} There are 101 coins, and only one of them differs from the other (genuine) ones by weight. We have to determine whether this counterfeit coin is heavier or lighter than a genuine coin. How can we do this using two weighings?
\end{puzzle}
\end{quote}

\textbf{Solution.} First, let us try this for other numbers $N$ of coins that are easier to think about than 101. If $N$ is divisible by 3, we divide coins into three piles with the same number of coins. Then we compare the first pile against the second in the first weighing, and the first pile against the third in the second weighing. If both of the weighings unbalance, the fake coin is in the first pile, and we know whether it is heavier or lighter. If one of the weighings balances, then the balancing piles do not contain the fake coin. Again, we know the answer. 

If $N$ is divisible by four, we divide the coins into four equal piles. We compare the first and the second pile against the third and the fourth. The weighing has to unbalance. We then compare the first pile against the second. If the second weighing balances, the coin must be in the third or fourth pile, and we know whether it is heavier of lighter. If the second weighing unbalances, the coin must be in the first or second piles, and again we know whether it is heavier or lighter.

For any number of coins we can try to merge the two solutions together. We divide all coins into three piles of size $a$, $a$ and $b$, where $a \leq b \leq 2a$. In the first weighing we compare the first two piles. If they balance, then the fake coin must be among the remaining $b$ coins. Now pick any $b$ coins from both pans in the first weighing and compare them to the remaining $b$ coins. If the first weighing is unbalanced, then the remaining coins have to be real. For the second weighing we pick $a$ coins from the remaining pile and compare them to one of the pans in the first weighing.

The solution I just described doesn't cover the cases of $N$ equal 1, 2, or 5. The problem can't be solved for one or two coins. For five coins, we can compare two coins against two coins in the first weighing. If it balances, the fake coin is the leftover coin and comparing it to any real coin provides the answer. If the first weighing unbalances, the fake coin is one of the four coins on the pan and we can proceed as we described above with the algorithm for the number of coins divisible by 4.

Back to logic. To continue an analogy, I have to invent a puzzle in which we do not need to find liars, but only need to say if the liars exaggerate or diminishe.

\begin{quote}
\begin{puzzle} On this island of knights and knaves the knaves are all of the same type: either exaggerating or diminishing. You meet five islanders. 
A: There are three truth-tellers among us.
B: There are three liars among us.

What type of liars live on this island?
\end{puzzle}
\end{quote}

\textbf{Solution.} The total number of people by their account is six. So one of A or B is a truth-teller, and the other is an exaggerating liar.

Going back to alternator coins: the alternator that switches between real and heavier is not interesting to study, as by symmetry it is the same as the alternator that switches between real and lighter. But we can have an alternator that switches between lighter and heavier. We discuss this case in the next section.

\section{Lighter-Heavier Coin}

Let us consider a fake coin that changes its weight from heavier than real to lighter than real and back. It never weighs the same as the real coin. We will call it an \textit{oscillating} fake coin. A logic analogue would be a liar who switches between exaggerating and diminishing. We call such a liar an \textit{oscillating liar}.

\begin{quote}
\begin{puzzle} On this island there are two types of people: knights and oscillating liars. You meet two islanders A and B. A says, ``There is one truth-teller among us.'' B says, ``There is one truth-teller among us.''

Who are A and B?
\end{puzzle}
\end{quote}

\textbf{Solution.} They can't both be truth-tellers. Therefore, at least one of them is an oscillating liar. And as they both made the same statement, the other one is also a liar.

\begin{quote}
\begin{puzzle} There are 12 identical-looking coins. One of them is an oscillating fake. Find the fake coin and identify whether it starts in the light or heavy state by using a balance scale three times. 
\end{puzzle}
\end{quote}

Why are there 12 coins? Does it remind you of anything? This puzzle can be solved with a wonderful argument that reduces its solution to the solution of Puzzle~\ref{puzzle:12coins}. Suppose we have a strategy and each time we have an imbalance we mentally reverse the direction of the imbalance. That means the coin doesn't change state. We just need to find the fake coin that is not known whether it is heavier or lighter. We can use exactly the same strategy as before. The argument works for any number of coins.  

If we do not care about the starting state the problem is equivalent to finding the fake coin in Puzzle~\ref{puzzle:13coins} when there is no need to identify whether it is heavier or lighter.

Now that we showed that the lighter-heavier coin problem is reducible to a known problem, we can invent a more complicated fake coin that changes its behavior according to a three-cycle, such as from lighter to heavier to real. Potentially, we can also study the lighter-heavier-real coin, but due to symmetry it is enough to study only one of the cases. 

This case is discussed in detail in our paper \cite{KK}: the adaptive strategy matches the bound provided by the count of the outcomes,  but the best oblivious strategy for the unknown starting state is unknown.

Back to logic. An analogue of the lighter-heavier-real coin is a person who makes statements with period three: a diminishing lie, an exaggerating lie and the truth. Let us call such a person a 3-rotator. It is easy to find the state of such a person by asking: How much is two plus two.

\begin{quote}
\begin{puzzle} On an island with truth-tellers and 3-rotators, you meet two islanders A and B. A says, ``There are zero truth-tellers between us.'' Who is who?
\end{puzzle}
\end{quote}

I leave it to the reader to solve this puzzle. In addition, the reader might try to identify the state of the 3-rotators if there are any.

\section{A Lying Scale}

We were discussing how the coins are like people: lying or telling the truth. In the coin-weighing problems, there is another player: the scale. The scale can also lie. 

But scales are not binary, so we have to be more specific to define a lying scale. Here is the simplest definition. The first type of the lying scale produces the opposite output of the truth-telling scale. It produces a balance when the coins on both pans weigh the same, but when there is an imbalance, the scale tilts the wrong way. We call such a scale---a \textit{reverse} scale. Let us call the scale that shows the correct weighings a \textit{true} scale.

\begin{quote}
\begin{puzzle} We have $N$ coins; one is fake and lighter and one scale that could be true or reverse. How many weighings do we need to decide whether the scale is true?
\end{puzzle}
\end{quote}

\textbf{Answer.} This is a no-brainer. We can put one coin on one pan and leave the other pan empty. The answer is one weighing.

Let's make the puzzle more mathematical. Suppose we only allow the same number of coins on both pans.

\begin{quote}
\begin{puzzle} We have $N$ coins; one is fake and lighter and one scale that could be true or reverse. How many weighings do we need to decide whether the scale is true if we are only allowed to put the same number of coins on the pans in one weighing?
\end{puzzle}
\end{quote}

\textbf{Solution.} The magic works again: there is a parallel between changing a fake coin from lighter to heavier and changing the scale from true to reverse. The algorithm is the same as in Puzzle~\ref{puzzle:h-or-l} where we are finding whether one fake coin is heavier or lighter. That means two weighings are enough if we have more than a total of two coins.

We can also have a different definition of lying. For example, we can invent a scale that always produces a wrong result, any wrong result. We call it a \textit{lying} scale. Again by leaving one pan empty and putting a coin on the other pan, we can always differentiate between a true scale and a lying scale. The next question is trickier. 

\begin{quote}
\begin{puzzle} We have $N$ coins, one of which is fake and lighter. We also have a scale that could be true or lying. Given that we only allow the same number of coins on both pans during a weighing, how many weighings do we need to differential between a true and lying scale?
\end{puzzle}
\end{quote}

I leave it to the readers to analyze this case.

\section{A Random Scale}

Let’s consider a scale that behaves like a normal person: sometimes lies, sometimes tells the truth. But we shouldn't call this scale a normal scale, as it is not normal at all. We call this scale a \textit{random} scale. The scale randomly decides the output of the weighing. A scale like that is a more realistic object than a coin that changes weight. Not surprisingly, the following puzzle appeared before chameleons were invented. I do not know the origins of this puzzle, but Knop \cite{KnopSlyScales} wrote a very good paper on this subject in Russian.

\begin{quote}\label{puzzle:randomscale}
\begin{puzzle} We have $3^{2n}$ identical-looking coins: one is fake and lighter. There are three balance scales, one of which is random. Given that the random scale is adversarial, what is the smallest number of weighings needed to find the fake coin?
\end{puzzle}
\end{quote}

\textbf{Solution idea.} We do not know which scale is playing against us. That means we can't avoid having at least one weighing that is lying. On the other hand we can try to set up a system that catches a lie as soon as possible. As soon as the random scale is identified we can use another scale to find the fake coin as we already know how to do. In all solutions I've seen, the worst case is when the random scale gives a correct answer until the very end, and then lies.

\textbf{The bound.} As before, we can draw a parallel between this problem and another problem that is much easier. Suppose the random scale lies not more than once. Then we can compare this problem to a problem where we have one scale that is allowed to lie not more than once. Now we are all set to produce a bound. Suppose there is a strategy in this case that requires $w$ weighings. Then the number of possible outcomes is $3^w$. After we have found the fake coin, we can look back at all the weighings and we will know which weighing lies. Any one of $w$ weighings can lie, and there are two possibilities for what the lying weighing can show. Plus it is possible that there are no lies. That means $2w+1$ different outcomes can point to the same coin. It follows that the number of coins we can process is not more than $3^w/(2w+1)$.

The paper \cite{KnopSlyScales} shows a strategy to find the fake coin in Puzzle~\ref{puzzle:randomscale} in $3n+1$ weighings. It might be possible to do better.

I actually heard about this puzzle with a particular value of 729 coins circulating at MOP (Math Olympiad Summer Program). The bound tells us that we need at least 9 weighings. The puzzle is to find a solution in 9 weighings. I leave the solution to the readers, but mention that our bound shows that 8 weighings is not enough.

Now that we have lying scales and lying coins, is there a parallel in logic that combines both? We can say that scales are people and coins are facts that they know. A fake coin represents a wrong fact. Scales describe their opinions about facts that might be wrong. Ideas like this were already introduced in logic. See the next section.

You might think that sane liar and an insane truth-teller can't be differentiated. This is why you have more puzzle tp think about.

\begin{quote}
\begin{puzzle} How will a sane liar and an insane truth-teller answer the question: ``Do you believe that two plus two is four?''
\end{puzzle}
\end{quote}

\section{Sane and Insane People}

Raymond Smullyan wrote a lot of books with logic puzzles. In his book \textit{The Lady or the Tiger?} \cite{RSLOrT} he introduced \textit{sane} and \textit{insane} people. 

Sane people believe that $2 + 2 =4$, and insane that $2+2 \neq 4$. In other words, sane people believe true statements and insane false statements.

A chapter of Smullyan’s book is devoted to an asylum, in which each member is either a doctor or a patient. All the doctors are supposed to be sane and all the patients insane. Inspector Craig is sent to inspect the asylum. 

\begin{quote}
\begin{puzzle} What statement can a sane patient make to prove that they should be released?
\end{puzzle}
\end{quote}

\textbf{Answer.} I am a patient.

You might assume that insane people are similar to liars. Here is another puzzle.

\begin{quote}
\begin{puzzle} A member of the asylum says, ``I believe I am a doctor.'' Can we conclude who they are?
\end{puzzle}
\end{quote}

\textbf{Answer.} It might seem that we can't. Indeed, both insane patients and sane doctors would say, ``I am a doctor.'' But this statement is different. An insane patient believes he is a doctor. That means, he believes that he doesn't believe he is a doctor. That means, an insane patient would say, ``I believe I am a patient.''

Now if we allow the members of the asylum to lie, we get the following. A sane liar and an insane truth-teller will both say that two plus two is not four. Compare this to the fact that a true scale produces an imbalance in the presence of a fake coin and a lying scale produces an imbalance when all the coins are real.

\section{Limitations of Scales}

People are more flexible than scales; they are able to talk about a lot of things. The scales only compare weights. That means there is a bigger variety of logic puzzles than coin-weighing puzzles. On the other hand, weighings are more structured and are easier to study.

In many logic puzzles people talk about each other. For example, see Puzzle~\ref{puzzle:kk}.

But scales can't talk about other scales. Or can they?

\begin{quote}
\begin{puzzle} There are three identical-looking balance scales that could be true or reversed. There are also $N$ coins and one of them is fake and lighter. What is the smallest number of weighings that you need to figure out the type of each scale?
\end{puzzle}
\end{quote}

\textbf{Solution.} Surprisingly, you can do it in one weighing if you put scales on scales. Here is the solution: place one coin on one of the pans of scale one; place two coins on one of the pans of scale two; place scales one and two on opposite pans of scale three. This way each weighing is an imbalance, and you know how a true scale should interpret it.


\begin{thebibliography}{9}

\bibitem{STEPAlternators}
B.~Chen, E.~Erives, L.~Fan, M.~Gerovitch, J.~Hsu, T.~Khovanova, N.~Malur, A.~Padaki, N.~Polina, W.~Sun, J.~Tan, and A.~The, Alternator Coins, \textit{Math Horizons}, v.25.1, pp. 22--25, (2017).

\bibitem{Dyson} F.~J.~Dyson, Note 1931---The problem of the pennies, Math. Gaz., 30 (1946) 231--234.

\bibitem{FGI} D.~Fomin, S.~Genkin, and I.~Inteberg, \textit{Mathematical Circles: Russian Experience},
 (Mathematical World, Vol. 7) 

\bibitem{TKRSblog} Tanya Khovanova, Raymond Smullyan's Magic Trick, Tanya Khovanova's Math Blog, available at \url{http://blog.tanyakhovanova.com/2010/05/raymond-smullyans-magic-trick/}

\bibitem{KK}
T.~Khovanova and K.~Knop, Coins that Change Their Weights, arXiv:1611.09201 [math.CO], (2016).

\bibitem{KKP}
T.~Khovanova, K.~Knop and O.~Polubasov, Chameleon Coins, arXiv:1512.07338 [math.HO], (2015).

\bibitem{KnopSlyScales}
K.~Knop, Weighings on ``sly scales'', (in Russian), Matematika v shkole, (2009) v2.

\bibitem{KP} K.~Knop, O.~Polubasov, Two counterfeit coins revisited, (In Russian) available at \url{http://knop.website/math/ff.pdf}, (2015)

\bibitem{OEIS}
\textit{The On-Line Encyclopedia of Integer Sequences}, published electronically at \url{https://oeis.org}, (2010), Sequence A001045.

\bibitem{Schell} E.~D.~Schell, Problem E651---Weighed and found wanting, Amer. Math. Monthly, 52 (1945) 42.

\bibitem{RSWhat} Raymond Smullyan, \textit{What is the Name of this Book?}, Prentice-Hall, (1978)

\bibitem{RSLOrT} Raymond Smullyan, \textit{The Lady or the Tiger?}, Knopf, (1982)

\end{thebibliography}
\end{document}